\documentclass{article}
\usepackage{arxiv}
\usepackage[unicode, pdftex]{hyperref}
\usepackage{url}             
\usepackage{booktabs}        
\usepackage{amsfonts}        
\usepackage{nicefrac}        
\usepackage{microtype}       
\usepackage{graphicx}
\usepackage{amsmath}
\usepackage{amssymb}  
\usepackage{natbib}
\usepackage{doi}
\usepackage{xcolor}
\usepackage[all]{xy}
\usepackage{caption}
\usepackage{subcaption} 
\usepackage{orcidlink}
\usepackage{tikz}
\usepackage{amsthm}
\usepackage{natbib}
\usepackage{algorithm}
\usepackage{algpseudocode}  
\usepackage{float}
\newcounter{cdn}
\tikzset{B/.style={blue!50!black},R/.style={red}}

\newtheorem{theorem}{Theorem}[section]

\theoremstyle{definition}
\newtheorem{definition}[theorem]{Definition}
\theoremstyle{remark}
\newtheorem{remark}[theorem]{Remark}
\newtheorem{example}[theorem]{Example}

\title{Complete Topological Invariant and Coding of Homoclinic Flows on the Two-Dimensional Sphere}

\renewcommand{\shorttitle}{Complete Topological Invariant and Coding of Homoclinic Flows}

\author{
  Alexandr Prishlyak \orcidlink{0000-0002-7164-807X}
  \\
  Taras Shevchenko National University of Kyiv, Ukraine
  \\
  \texttt{prishlyak@knu.ua}
  \And
  Illia Ovtsynov \orcidlink{0009-0000-3262-2291} 
	\\
  Taras Shevchenko National University of Kyiv, Ukraine \\
  \texttt{iliarkov@gmail.com}
}

\hypersetup{
pdftitle={Complete Topological Invariant and Coding of Homoclinic Flows on the Two-Dimensional Sphere},
pdfauthor={Alexandr Prishlyak, Illia Ovtsynov},
pdfkeywords={Morse flow, topological classification, vector field, signed chord diagram, dual diagram, enumeration algorithm, surface topology},
colorlinks=true,
linkcolor=blue,
citecolor=blue,
urlcolor=blue
}
\begin{document}
\maketitle

\begin{abstract}
%
We study codimension one homoclinic flows on closed surfaces and introduce a complete topological invariant based on the structural graph of the flow and the distinguishing graphs of its atoms. This invariant completely determines the topological equivalence class of a homoclinic flow. For flows on the two-dimensional sphere, we construct a canonical coding algorithm that avoids the global renumbering required in the classical coding of Morse flows and substantially reduces the computational complexity. As an application, we derive recursive formulas for the number of topologically distinct homoclinic flows on the two-dimensional sphere and compute these numbers for flows with a small number of saddle singularities and closed orbits.
\end{abstract}






\section{Introduction}
\label{sec:1}

Morse--Smale vector fields form the class of structurally stable dynamical systems and play a crucial role in the qualitative theory of dynamical systems. In 1970, Palis and Smale \cite{Palis1970} proved that Morse--Smale vector fields form an open and dense subset in the space of vector fields on closed surfaces.

The simplest structurally unstable vector fields are the codimension one vector fields, which appear in generic one-parameter families of vector fields. They form the codimension one strata in the space of vector fields and belong to the boundary of the Morse--Smale vector fields set. 

In this paper, we study codimension one flows on closed surfaces possessing exactly one homoclinic saddle loop. We are interested in their classification up to topological trajectory equivalence, which will simply be referred to as \emph{topological equivalence}. Since merely the orbit structure is relevant, and we are considering smooth vector fields, we shall freely identify vector fields with the flows generated by them and use the terms \emph{vector field} and \emph{flow} interchangeably throughout the paper.

Complete topological invariants are important in the topological classification of Morse--Smale flows and their codimension one bifurcations. Such invariants are usually described by graphs equipped with additional combinatorial data. Classical examples include the separatrix diagram of Peixoto \cite{Peixoto1973}, embedded graphs encoded by rotation systems \cite{GT87}, and chord diagrams for Morse flows with one source \cite{Kybalko2018}. This problem was also thoroughly studied in \cite{Oshemkov1998}. Flows with collective dynamics were researched from a similar point of view in \cite{prus21}.
Related problems concerning the topology of functions on surfaces and their generic deformations were studied in \cite{Hladysh2019}, where deformations in the general position of optimal functions on oriented surfaces with boundary were investigated.

In our recent paper \cite{provts}, we constructed a complete topological invariant for gradient codimension one flows on closed surfaces, and in \cite{provts26} for compact surfaces with boundary. Those flows contain neither periodic orbits nor homoclinic loops. The invariant extends the chord diagram of Morse flows by introducing a distinguished point, which provides a canonical numbering of the chords and, consequently, a canonical code of the flow.

The aim of the present paper is to construct a complete topological invariant for codimension one homoclinic flows on closed surfaces. In contrast to the gradient case, such flows possess a unique homoclinic saddle loop and may contain periodic orbits. These flows were researched in \cite{Homburg_2013} and \cite{Kibkalo_2022}. The proposed invariant is again based on a chord diagram equipped with additional combinatorial data adapted to the homoclinic loop. As an application, we derive explicit formulas for counting topologically distinct codimension one homoclinic flows on the sphere.

Throughout this paper, a \emph{graph} is understood as a finite one-dimensional CW-complex (topological graph). In particular, loops and multiple edges are allowed.

\section{Homoclinic Flows on Closed Surfaces}

Throughout this paper we are assuming that $M$ is the two-dimensional sphere and $$g\colon M \times \mathbb{R} \to M$$ is a $C^{\infty}$ flow on a closed surface $M$. Furthermore, we denote $\Omega_0(g)$ a set of all the sources of $g$, $\Omega_1(g)$ a set of all the saddles of $g$ and $\Omega_2(g)$ a set of all the sinks of $g$.

\begin{definition}
  A set $$W^+(p) = \big\{x \in M\,|\,\lim\limits_{t \to +\infty} g(x, t) = p\big\}$$ is called a \emph{stable manifold} of a singular point $p$ of the flow $g$ whereas a set $$W^{-}(p) = \big\{x \in M\,|\,\lim\limits_{t \to -\infty} g(x, t) = p\big\}$$ is called an \emph{unstable manifold} of a singular point $p$.
\end{definition}

\begin{definition}
Let $u_0$ be a saddle point of $g$. If there exists a trajectory $\gamma$ of $g$ such that
\[
\alpha(\gamma)=\omega(\gamma)=u_0,
\]
then $u_0$ is called a \emph{homoclinic saddle}, and $\gamma$ is called a \emph{homoclinic trajectory}.
\end{definition}

We also state a neighbourhood $U$ of a subset $S$ of $M$ to be \emph{regular} in case the singular points of $g$ belonging to $S$ are the only singular points belonging to $U$ and $\partial U$ does not intersect any trajectory of $g$ more than twice.

Let $\gamma$ be a homoclinic trajectory of a saddle $u_0$. Put $L=\gamma\cup\{u_0\}$.

The closure of the component containing trajectories for which $L$ is either the $\alpha$-limit set or the $\omega$-limit set is called the \emph{$H$-component} (where $H$ stands for \emph{hyperbolic}, since this component contains one hyperbolic sector). The closure of the other component is called the \emph{$DS$-component} (standing for \emph{double saddle}, since it is a two-dimensional disk which is topologically equivalent to the one with a double saddle on its boundary).

In this paper, we study codimension-one homoclinic flows on the boundary of Morse flows on the two-dimensional sphere. Such flows naturally arise in one-parameter bifurcations of Morse vector fields and can be approximated by Morse vector fields. We introduce a complete topological invariant for these flows, define an equivalence relation on these invariants, and construct a finite code which determines the invariant, and hence the flow itself, up to topological equivalence.

\begin{definition}
A $C^{\infty}$ vector field on $M$ is called an \emph{$HC_{+}$-vector field} (respectively, an \emph{$HC_{-}$-vector field}) if it satisfies the following conditions:

\begin{enumerate}
\item[(1)] It has finitely many singular points and periodic orbits all of which are hyperbolic.


\item[(2)] For every trajectory contained in the $DS$-component, both its $\alpha$- and $\omega$-limit sets are critical elements.

\item[(3)] For every trajectory contained in the interior of the $H$-component, its $\omega$-limit set (respectively, $\alpha$-limit set) is a critical element, whereas its $\alpha$-limit set (respectively, $\omega$-limit set) is either a critical element or the boundary of the $H$-component.

\item[(4)] The vector field contains a homoclinic trajectory $\gamma$.

\item[(5)] If
\[
x\in\left(W^{u}(u_1)\cap W^{s}(u_2)\right)\setminus\gamma,
\]
where $u_1,u_2\in\operatorname{Sing}(g)$, then
\[
T_xW^{u}(u_1)+T_xW^{s}(u_2)=T_xM.
\]
\end{enumerate}
\end{definition}

\begin{definition}
A $C^{\infty}$ flow on a closed surface $M$ is called an \emph{$HC_{+}$-flow} (respectively, an \emph{$HC_{-}$-flow}) if it is generated by an $HC_{+}$-vector field (respectively, an $HC_{-}$-vector field).
\end{definition}

$HC_{+}$-flows and $HC_{-}$-flows will be referred to collectively as
\emph{homoclinic flows}. These are precisely the codimension one flows
considered throughout this paper.

\section{Complete Topological Invariant}


\subsection{Atoms and Distinguishing Graphs}

Let $$g\colon M \times \mathbb{R} \to M$$ be an $HC_{-}$-flow on $M$ with a homoclinic saddle $u_{0}$ and homoclinic trajectory $\gamma$.

Denote by $A^{-}$ the disjoint union of arbitrarily fixed regular neighborhoods of all sources and repelling periodic orbits of $g$, and by $A^{+}$ the disjoint union of arbitrarily fixed regular neighborhoods of all sinks and attracting periodic orbits of $g$. Let $U_{S}$ be a fixed regular neighborhood of the set
$S=\gamma\cup\{u_{0}\}$.

The periodic orbits of $g$ decompose the surface $M$ into connected regions. The closures of these regions are called \emph{atoms}. The unique atom containing the set $S$ is called the \emph{$H$-atom}, while all remaining atoms are called \emph{$C$-atoms}.

The homoclinic trajectory $\gamma$ divides the $H$-atom into two connected components. The closure of the component whose interior is disjoint from $W^{+}(u_{0})\cup W^{-}(u_{0})$ is called the \emph{$H_{-}$-atom}; the closure of the other component is called the \emph{$H_{+}$-atom}.

Henceforth, if we speak of \emph{atoms}, we will mean only $C$-atoms, $H_+$-atoms, or $H_-$-atoms. Moreover, we agree that an orientation on each atom is defined by the direction of the homoclinic trajectory.

Also it is worth drawing one's attention that we will define for a specific subset $A \subset M$ a flow $$g|_{A \times \mathbb{R}}\colon A \times \mathbb{R} \to A$$ so that for any $t \in \mathbb{R}$, $x \in A$ we obtain $g|_{A \times \mathbb{R}}(x, t) = g(x, t)$ provided the equation does not arise any contradictions.

We first define the distinguishing graph of a $C$-atom.

\begin{definition}
Let $L\subset M$ be a $C$-atom. The \emph{distinguishing graph} of $L$ is the disconnected union
\[
\Gamma_{L}=\Gamma_{+}\sqcup\Gamma_{-},
\]
where
\[
\Gamma_{\pm}
=
(\partial A^{\mp}\cap L)
\cup
\bigcup_{s\in\Omega_{1}(g|_{L\times\mathbb{R}})}
\left(
W^{\pm}(s)\cap(L\setminus A^{\mp})
\right),
\]
and the following conditions hold:

\begin{enumerate}
\item[(1)]
The vertices of $\Gamma_{\pm}$ are the points of
\[
V_{\pm}
=
\partial A^{\mp}
\cap
\bigcup_{s\in\Omega_{1}(g|_{L\times\mathbb{R}})}
W^{\pm}(s),
\]
while the edges are the connected components obtained after removing the vertices.

\item[(2)]
The vertices of $\Gamma_{+}$ are colored blue, whereas the vertices of $\Gamma_{-}$ are colored red.

\item[(3)]
The edges from $\partial A^+ \cap L$ ($\partial A^- \cap L$) are colored thick red (thick blue) whereas the other edges in $\Gamma_-$ ($\Gamma_+$) are colored thin red (thin blue).

\item[(4)]
Every edge of the distinguishing graph lying on the boundary of a neighborhood of a periodic orbit is oriented according to the direction of the flow.
\end{enumerate}
\end{definition}

The connected components of the union of all oriented edges of a distinguishing graph are called its \emph{oriented cycles}. All the connected components of the union of all the thick red and thick blue edges of a distinguishing graph are called its \emph{boundary cycles}.

\begin{definition}
  Let $L \subset M$ be an $H_+$-atom (or an $H_-$-atom). The \emph{distinguishing graph} of $L$ is a union $\Gamma_L$ of a disconnected union of graphs $\Gamma_+$ and $\Gamma_-$ and the set $U := \partial S \cap L$ where
  \[
  \Gamma_{\pm} = (\partial A^{\mp} \cap L) \cup \bigcup\limits_{s \in \Omega_1(g|_{L \times \mathbb{R}})} (W^{\pm}(s) \cap (L \setminus A^{\mp}) \cap (L \setminus S)),
  \]
  and the following conditions hold:

\begin{enumerate}
\item[(1)]
The vertices of $\Gamma_{\pm}$ are the points of
\[
V_{\pm} = (\partial A^{\mp} \cup U) \cap \bigcup\limits_{s \in \Omega_1(g|_{L \times \mathbb{R}})} W^{\pm}(s),
\]
while the edges are the connected components obtained after removing the vertices.

\item[(2)]
The vertices of $\Gamma_{+}$ are colored blue, whereas the vertices of $\Gamma_{-}$ are colored red, and the set $U$ is colored black.

\item[(3)]
The edges from $\partial A^+ \cap L$ ($\partial A^- \cap L$) are colored thick red (thick blue) whereas the other edges in $\Gamma_-$ ($\Gamma_+$) are colored thin red (thin blue).

\item[(4)]
Every edge of the distinguishing graph lying on the boundary of a neighborhood of a periodic orbit or a homoclinic trajectory is oriented according to the direction of the flow.
\end{enumerate}
\end{definition}

As in the definition before, the connected components of the union of all oriented edges of a distinguishing graph are called its \emph{oriented cycles}. All the connected components of the union of all the thick red, thick blue and black edges of a distinguishing graph are called its \emph{boundary cycles}.

We are now ready to define the distinguishing graph of the entire flow.

\begin{definition}
The \emph{distinguishing graph} of an $HC_{-}$-flow $g$ on a closed surface $M$
is the collection consisting of the distinguishing graphs of the
$H_{+}$-atom, the $H_{-}$-atom, and all $C$-atoms of $g$; besides, following conditions hold true:
\begin{enumerate}
  \item[(1)] The non-black oriented cycles occurring in these distinguishing graphs are labeled by
the integers $1,\ldots,N$, where $N$ is the number of periodic orbits of $g$ whereas the black oriented cycle is label by $0$.
  \item[(2)] Two non-black oriented cycles receive the same label if and only if they belong to
distinguishing graphs of different atoms and correspond to the same periodic
orbit of the flow.
\end{enumerate}
\end{definition}

\begin{definition}
Two finite collections of graphs are said to be \emph{isomorphic} if they
contain the same number of graphs and there exists a bijection between the two
collections such that corresponding graphs are isomorphic.
\end{definition}

\begin{theorem}\label{main}
Let $g_{1}$ and $g_{2}$ be two $HC_{-}$-flows on $M$, each
having $N$ periodic orbits. Then $g_{1}$ and $g_{2}$ are topologically
equivalent if and only if their distinguishing graphs
$\Lambda_{1}$ and $\Lambda_{2}$ satisfy the following conditions.

\begin{enumerate}
\item[(1)]
The collections $\Lambda_{1}$ and $\Lambda_{2}$ are isomorphic, and the
corresponding graph isomorphisms preserve the colors of the edges and the
orientations of the oriented cycles.

\item[(2)]
There exists a permutation
\[
\sigma:\{0,\ldots,N\}\rightarrow\{0,\ldots,N\}
\]
with its rest point $0$ such that if an oriented cycle of a graph in $\Lambda_{1}$ is labeled by
$i \in \{1,\ldots,N\}$, then the corresponding oriented cycle of the corresponding graph in
$\Lambda_{2}$ is labeled by $\sigma(i)$.
\end{enumerate}
\end{theorem}

\subsection{Structural Graph}


Now let us construct a structural graph of an $HC_{-}$-flow $g$ on the $2$-dimensional sphere $M$.

Suppose that $g$ has $N$ $C$-atoms, $N \in \{0\} \cup \mathbb{N}$. We denote them as $C_1, \ldots, C_N$. Also we have an $H_+$-atom $H_+$ and an $H_-$-atom $H_-$.

We fix in an arbitrary way points $x_i \in C_i$, $i = \overline{1,N}$, and $y_+ \in H_+$, $y_- \in H_-$.

\begin{definition}
  A \emph{structural graph} of the flow $g$ is a topological graph for which following conditions hold true:
  \begin{enumerate}
    \item[(1)] $x_1, \ldots, x_N, y_+, y_-$ are the vertices of this graph and there are a symbol `+' near $y_+$ and a symbol `-' near $y_-$;
    \item[(2)] two vertices of this graph are connected with a single edge if and only if the corresponding atoms on $M$ are adjacent;
    \item[(3)] if a vertex belongs to an $C$-atom, then it is colored blue, otherwise it is colored black;
    \item[(4)] if an edge connects two black vertices it is colored black and called an \emph{H-edge};
    \item[(5)] if an edge is not an H-edge then it is colored red if it connects two vertices corresponding to atoms divided by an attracting cycle and it is colored blue if it connects two vertices corresponding to atoms divided by a repelling cycle;
    \item[(6)] edges which are not black are assigned with numbers from $1$ to $N$ in a random way;
    \item[(7)] there is a symbol $``+``$ or $``-``$ near each edge provided that $``+``$ is placed iff the appropriate closed or homoclinic trajectory is directed along the defined orientation.
  \end{enumerate}
\end{definition}

A structural graph is actually a \emph{tree} --- a connected graph without cycles.






\begin{theorem}
  Let $\Lambda$ be a distinguishing graph for an $HC_-$-flow $g$ on $M$. Let also $\Sigma$ be a tree in which vertices and edges are colored either black or blue. Then $\Sigma$ is a structural graph for $g$ if there is a bijection $f$ mapping a set of vertices in $\Sigma$ into a set of distinguishing graphs in $\Lambda$ and there is a permutation $\sigma:\{1,\ldots,N\}\rightarrow\{1,\ldots,N\}$ for which next conditions hold true:
  \begin{enumerate}
    \item[(1)] blue vertices in $\Sigma$ are mapped into graphs for $C$-atoms for $\Lambda$ whereas the black vertex with the symbol $``+``$ is mapped into the graph for the $H_+$-atom and the black vertex with the symbol $``-``$ is mapped into the graph for the $H_-$-atom;
    \item[(2)] two vertices in $\Sigma$ are connected with a black edge if and only if corresponding graphs in $\Lambda$ have black edges;
    \item[(3)] two vertices in $\Sigma$ are connected with a red edge labeled by number $i \in \{1,\ldots,N\}$ if and only if corresponding graphs in $\Lambda$ have red oriented cycles having the same label $\sigma(i)$;
    \item[(4)] two vertices in $\Sigma$ are connected with a blue edge labeled by number $i \in \{1,\ldots,N\}$ if and only if corresponding graphs in $\Lambda$ have blue oriented cycles having the same label $\sigma(i)$;
    \item[(5)] there is always $``+``$ near the black edge in $\Sigma$ and any other edge labeled by number $i \in \{1, \ldots, N\}$ is marked by $``+``$ iff the corresponding graphs have oriented cycles having the same label $\sigma(i)$ which assign the direction of the appropriate closed trajectories with respect to the defined orientation on $M$.
  \end{enumerate}
\end{theorem}


But it is important to highlight that actually Theorem~\ref{main} provides a full list of conditions which help identify whether two $HC_-$-flows are topologically equivalent or not. Furthermore, having a distinguishing graph for an $HC_-$-flow one has an opportunity to recover the flow according to its topological nature.

Nevertheless, it is not enough to recover completely such a flow having only its structural graph at disposal --- one can only understand how periodic orbits of the flow are placed on the sphere, the behaviour of trajectories near them and which atom contains a homoclinic trajectory. It can be helpful in some cases but not when the overall topological characteristics of the flow is needed. Apart from that, a structural graph can be obtained from the distinguishing graph (but not on the contrary).

\section{Encoding of Homoclinic Flows}

Let $$g\colon M \times \mathbb{R} \to M$$ be an $HC_-$-flow with $N$ closed trajectories, $N \in \{0\} \cup \mathbb{N}$, on the two-dimensional sphere $M$.

After the speculations from the previous section, we now can easily construct a distinguishing graph $\Lambda$ for this flow having $N$ oriented cycles. Furthermore, we are able to construct a structural graph $\Sigma$ for $g$.

For convenience we will make sure that an edge in a structural graph and a pair of oriented cycles which are corresponding to the same closed trajectory in $g$ are assigned with the same number $i \in \{1, \ldots, N\}$. It appears that now each $C$-atom is assigned unambiguously with a specific natural number.

Let $B$ be the set of blue vertices and $C$ be the set of black vertices in $\Sigma$. We are about to construct a function for numbering the blue vertices in $\Sigma$ $$f\colon B \cup C \to \{``0+``, ``0-``, ``1``, ``2``, \ldots\}.$$

We will require a helpful notation for making our speculations plainer. Let $A$ be a set of vertices in a tree $G$ which have degree $1$. Denote $G_A$ the subgraph of $G$ which can be obtained through the exclusion of the vertices from $A$ and the edges incident on them.

Returning to $\Sigma$, we define $f$ according to the following algorithm:
\begin{enumerate}
  \item[(1)] on the first step, we choose a non-black edge in $\Sigma$ with an endpoint of degree $1$ which has the least assigned natural number $i_1$ (we denote it as $v_1$) --- thus, we state that $f(v_1) = ``i_1``$; the corresponding edge has colour $G_1$ (red or blue);
  \item[(2)] on the $n$-th step, we choose a non-black edge in $((\Sigma_{\{v_1\}})_{\{v_2\}})\ldots_{\{v_{n-1}\}}$ with an endpoint of degree $1$ which has the least assigned natural number $i_n$ (we denote it as $v_n$) --- thus, we state that $f(v_n) = ``i_n``$; the corresponding edge has colour $G_n$ (red or blue);
  \item[(3)] the procedure is done until all the blue vertices have been numbered --- afterwards we state that $f(x) = ``0+``$ and $f(y) = ``0-``$ provided that $x$ is a vertex belonging to the $H_+$-atom and $y$ is a vertex belonging to the $H_-$-atom.
\end{enumerate}



We are about to define a structural part for the flow code. Its idea is directly derived from the concept of construction of a Pr$\ddot{\textrm{u}}$fer sequence for standard trees but it is refined so it can be used for coding $HC_-$-flows in a proper way.

Now we construct two specific functions $$G\colon \{1, \ldots, N\} \to \{``0+``, ``0-``, ``1``, ``2``, \ldots\},\,\,I\colon \{1, \ldots, N\} \to \{``+``, ``-``\}.$$ To describe it less formally, our structural part of the code will have $N$ symbols in itself. These mappings will explain which exactly symbol will be placed in each position and which indices it will possess. So, for $i \in \{1, \ldots, N\}$ we claim that $$G(i) = f(w_i),\,\,I(i) = f$$ where $w_i$ is the other endpoint of the edge in the graph $((\Sigma_{\{v_1\}})_{\{v_2\}})\ldots_{\{v_{i-1}\}}$ incident to $v_i$; furthermore, the same edge in $\Sigma$ is assigned with the label $f \in \{``+``, ``-``\}$.

\begin{definition}
  A \emph{structural flow code} for the flow $g$ is a sequence of symbols $$S(\ldots)$$ where instead of dots the symbols $G(1), G(2), \ldots, G(N)$ assigned respectively with their lower right indices $I(1), I(2), \ldots, I(N)$ are consecutively written respectively in colours $G_1, G_2, \ldots, G_N$.
\end{definition}



\begin{remark}
If the flow has no $C$-atoms, then the structural flow code has no symbols inside round brackets.
\end{remark}


What is left to do is to write a code for distinguishing graphs for all the atoms. But we must memorize all the objects and notations made in this chapter before.




Let's fix an arbitrary atom $L$. We construct a flow code for its distinguishing graph $\Gamma$. 

We denote the set of all the vertices of $\Gamma$ as $T$ and the set of all the oriented cycles as $C$. Now, let $$\sigma_L\colon C \to \{``0+``, ``0-``, ``1``, ``2``, \ldots\}$$ be a mapping which maps each cycle from $C$ into $f(y)$ where $y$ is a vertex in $\Sigma$ corresponding to the atom adjacent to $L$ through exactly this cycle. Additionally, $\Lambda := \Lambda_1 \cup \ldots \cup \Lambda_4$ where
\begin{center}
  $\Lambda_1 = \{``\{``, ``\}_{0+}``, ``\}_{0-}``, ``\}_{1}``, ``\}_{2}``, ``\}_{3}``, \ldots\},\,\,\,\,\Lambda_2 = \{``[``, ``]``, ``\langle``, ``\rangle``\},$ \\
  $\Lambda_3 = \{``(``, ``)_{0+}``, ``)_{0-}``, ``)_1``, ``)_2``, ``)_3``, \ldots\},\,\,\,\,\Lambda_4 = \{``1``, ``2``, ``3``, \ldots\}.$ \\
\end{center}


We also fix a permutation $\xi: T \to T$ which maps an endpoint of a thin non-black edge in $\Gamma$ to the other its endpoint.

Also note that, of course, it may appear that $T = \varnothing$. We state at once that the flow code for this kind of atoms is an empty string. Hereafter, we assume $\Gamma$ has a non-empty set of vertices.



\begin{definition}
  A sequence of symbols $$q_1 q_2 \ldots q_m$$ from $\Lambda$ is called a \emph{an $HC$-like code list} if any two-digit number in it is separated by two black commas from both left and right side; furthermore, $q_1, q_m \in \Lambda_i$ where $1 \leq i \leq 2$ but $q_2, q_3, \ldots, q_{m-1} \in \Lambda_4$ so that $q_1$ contains an opening bracket and $q_m$ contains a closing bracket.
\end{definition}


\begin{definition}
  A sequence of symbols $$q_0 L_1 L_2 \ldots L_m q_1$$ where $m \in \mathbb{N}$, is called an \emph{HC-like code} if each $L_i$, $i = \overline{1,m}$, is an HC-like code list and the first inclusions of natural numbers in the code form an increasing sequence; moreover, $q_0, q_1 \in \Lambda_3$ so that $q_0$ contains an opening bracket and $q_1$ contains a closing bracket.
\end{definition}

\begin{definition}
  A \emph{cycle list} in an HC-like code $L_1 L_2 \ldots L_m$ is a sequence of symbols $$a_1 a_2 \ldots a_{n-1}$$ where $a_1, a_2, \ldots, a_{n-1} \in \Lambda_4$, $a_1$ and $L_{q_1}$, $a_1 \in L_{q_1}$, $q_1 = \overline{1,m}$, are fixed provided that:
  \begin{enumerate}
    \item[(1)] each $a_i$, $i = \overline{1,n}$, is a number considered with memorizing in which code list it is placed and where namely in this code list it is placed;
    \item[(2)] if $a_i \in L_{q_i}$, $i = \overline{1,n-1}$, then the next number $x$ in the list (if $a_i$ is the last one $x$ must be the first one) has its second location in the code --- the corresponding code list is $L_{q_{i+1}}$ and $a_{i+1}=x$;
    \item[(3)] $a_1 = a_n$ but for any $i = \overline{2,n-1}$ $a_1 \neq a_i$;
    \item[(4)] numbers among $a_1, a_2, \ldots, a_{n-1}$ having more than one digit in this sequence are separated by commas from both left and right sides.
  \end{enumerate}
\end{definition}

A set of cycle lists being a set of all cyclic permutations of the same cycle list is called a \emph{cycle}.

\begin{definition}
  A sequence of symbols $$q_0 L q_1$$ is called a \emph{cycle-like code list} for an $HC$-like code if $L$ is a cycle list for this code and $q_0, q_1 \in \Lambda_1$ so that $q_0$ contains an opening bracket and $q_1$ contains a closing bracket.
\end{definition}

\begin{definition}
  A sequence of symbols $$q_0 L_1 L_2 \ldots L_m q_1$$ is called a \emph{cycle-like code} for an $HC$-like code if $L_1, \ldots, L_m$ is list of some mutually distinct cycle-like code lists for the $HC$-like code and $q_0, q_1 \in \Lambda_3$ so that $q_0$ contains an opening bracket and $q_1$ contains a closing bracket.
\end{definition}

Now, we consider a subgraph $\Gamma_1$ of $\Gamma$ which contains all non-red edges and blue vertices.

Suppose that $M_1, \dots, M_n$ are boundary cycles of $\Gamma_1$. Let also $v_i^1 \in M_i \cap T$, $i = \overline{1, n}$, be some arbitrarily fixed point --- we call it a \emph{distinguishing vertex}. We have an ordered set $$V_i = \{v_i^1, v_i^2, \dots, v_i^{n_i}\} = T \cap M_i$$ obtained after bypassing $M_i$ completely along the fixed orientation starting from the point $v_i^1$.


We start bypassing each boundary cycle $M_i$ according to the defined motion, gradually increasing $i$ by $1$. 
Once we meet any point $v \in V_i$, $i = \overline{1,n}$, we define  in this point a value of a function $$h_{i, v_i^1}: V_i \to \mathbb{N} \cup \{0\},\,\,i = \overline{1,n},$$ which tells us endpoints of how many thin blue edges we have met since the beginning of bypassing $M_1$. Then we claim a value of a corresponding function $$f_{i, v_i^1}: V_i \to \Lambda_5,\,\,i = \overline{1,n},$$ in the same point. $f_{i, v_i^1}(v_i^n)$ corresponds to $m$, $n = \overline{1, n_i}$, $i = \overline{1,n}$, if:
\begin{enumerate}
     \item[(1)] $\xi(v_i^n) \notin M_1 \cup M_2 \cup \dots \cup M_{i-1} \cup (M_i \cap \{v_i^1, v_i^2, \dots, v_i^{n-1}\})$ and $h_{i, v_i^1}(v_i^n) = m$;
     \item[(2)] $\xi(v_i^n) = v_j^l$, $1 \leq j < i$, $1 \leq l \leq n_j$ where $f_{i, v_j^1}(v_j^l)$ corresponds to $m$;
     \item[(3)] $\xi(v_i^n) = v_i^q$, $1 \leq q < n$ where $f_{i, v_i^1}(v_i^q)$ corresponds to $m$.
\end{enumerate}

\begin{definition}\label{blue}
  Let $M_i$, $i = \overline{1,n}$, be a boundary cycle in $\Gamma_1$ and $v_i^1 \in M_i$ be its distinguishing vertex. A sequence of symbols from $\Lambda$
  $$q_1\,f_{i, v_i^1}(v_i^1)\,f_{i, v_i^1}(v_i^2)\,\ldots\,f_{i, v_i^1}(v_i^{n_i})\,q_2$$
  colored blue is called a \emph{code list for $M_i$ with a distinguishing vertex $v_i^1$} provided that:
  \begin{enumerate}
    \item[(1)] if $M_i$ is a thick blue edge which is not an oriented cycle then $q_1 = ``[``$ and $q_2 = ``]``$;
    \item[(2)] if $M_i$ is a thick blue oriented cycle then $q_1 = ``\{``$ and $q_2 \in \Lambda_1$ --- there is a string in the right lower corner of the closing curly bracket in $q_2$ which denotes $\sigma_L(M_i)$;
    \item[(3)] if $M_i$ is a black oriented cycle then $q_1 = ``\langle``$ and $q_2 = ``\rangle``$;
    \item[(4)] numbers having more than one digit in this sequence are separated by commas from both left and right sides.
  \end{enumerate}
\end{definition}


\begin{definition}\label{atom}
  Let $x_1 \in M_1$ be an arbitrary fixed vertex of $\Gamma_1$. A sequence of symbols $$\mu_1\,\lambda_1\,\lambda_2\,\ldots\,\lambda_n\,\mu_2$$ is called a \emph{flow code for the graph $\Gamma_1$ with a distinguishing vertex $x_1$} if each $\lambda_i$, $i = \overline{1,n}$, is a code list for a boundary cycle $M'_i$ with a distinguished vertex $v'_i \in M'_i$ for the flow $g$ so that:
  \begin{enumerate}
    \item[(1)] $M'_1 = M_1$, $v'_1 = x_1$;
    \item[(2)] if we take the least non-zero number in $\lambda_1\,\dots\,\lambda_i$ for each $i = \overline{1, n-1}$ included in it only once, it will represent some $v \in T$ --- then $v'_{i+1} = \xi(v)$ and $M'_{i+1}$ is a boundary cycle containing $\xi(v)$;
    \item[(3)] $\mu_1 = ``(``$ and $\mu_2 \in \Lambda_3$ ---  there is a string in the right lower corner of the closing round bracket in $\mu_2$ that defines $f(x)$ where $x$ is the vertex in $\Sigma$ corresponding to the atom $L$.
  \end{enumerate}

\end{definition}

Code lists for boundary cycles of $\Gamma_1$ with distinguishing vertices are actually $HC$-like code lists while the flow code for $\Gamma_1$ with a distinguishing vertex is actually an $HC$-like code.

After collapsing all the boundary cycles in the interior of the atom into points we will get that all the thin blue edges in $\Gamma_1$ divide it into different $2$-cells which we will call \emph{canonical $2$-cells}. Cycles, in fact, give us the boundaries of these $2$-cells in atoms --- a closure of each of them contains exactly one sink or one attracting closed trajectory or an attracting homoclinic trajectory if we speak of the $H_-$-atom.

To be more precise, there are several cycle lists which may correspond to a single canonical $2$-cell. So, we agree that we call a \emph{canonical cycle list} such a cycle list for such a $2$-cell which starts from the least possible natural number.

\begin{definition}
A sequence of symbols $$q_0 L q_1$$ colored red is called an \emph{advanced cycle} if $L$ is a canonical cycle list for a specific canonical $2$-cell with a red oriented cycle $R$ inside it, $q_0 = ``\{``$ and $q_1 \in \Lambda_1$ --- there is a string in the right lower corner of the closing curly bracket in $q_1$ which denotes $\sigma_L(R)$.
\end{definition}

The constructed code unambiguously defines a placement of stable manifolds of saddles inside of atoms. Moreover, during constructing the code we bypass different fragments of boundaries of these $2$-cells --- we will say in this case that we \emph{meet different cycles}.

The only thing left to define by the code is where we have exactly closed trajectories and how the orientation of these trajectories behaves. Thus, we have just to write down specific codes for the appropriate cycles. 

We define a linear order on the set of cycles in the atom $L$ --- for instance, a cycle $C_1$ is considered to be not greater than a cycle $C_2$ provided that the first meeting of a number from which the cycle $C_1$ is formed happens before the first meeting of a number from which the cycle $C_2$ is formed. 

\begin{definition}
  Suppose that the atom $L$ has $m$ different cycles, $m \in \mathbb{N}$. A sequence of symbols $$q_0 M L_1 L_2 \ldots L_m q_1$$ is called a \emph{cycle code} for the atom $L$ if:
  \begin{enumerate}
    \item[(1)] unless $L$ is the $H_-$-atom, $M$ is an empty string; otherwise it is an advanced cycle of $L$ corresponding to the homoclinic trajectory in its closure;
    \item[(2)] $L_1, \ldots, L_m$ are all its advanced cycles corresponding to all the attracting closed trajectories belonging to the closure of $L$ arranged in ascending order;
    \item[(3)] $\mu_1 = ``(``$ and $\mu_2 \in \Lambda_3$ ---  there is a string in the right lower corner of the closing round bracket in $\mu_2$ that defines $f(x)$ where $x$ is the vertex in $\Sigma$ corresponding to the atom $L$.
  \end{enumerate}
\end{definition}

We define a linear order on the set $\{``0+``, ``0-``, ``1``, \ldots, ``N``\}$ like that:

$$``0+`` \leq ``0-`` \leq ``1`` \leq ``2`` \leq \ldots \leq ``N``.$$

Now we are ready to write down two sections of a \emph{distinguishing flow code} for $g$ according to the further definition.

\begin{definition}
  A \emph{distinguishing flow code} for $g$ is two strings of symbols $$D_-\ldots$$ $$D_+\ldots$$ where the all non-empty flow codes for all the atoms containing blue symbols must replace the dots in the first string arranged in ascending order of their indices, and all the cycle codes for all the atoms of $g$ colored red and arranged in ascending order of their indices must replace the dots in the second string.
\end{definition}

We also call a distinguishing flow code is \emph{canonical} if the numbers in its first string (not the indices) are placed in a lexicographic order according to the canonical order of natural numbers.

\begin{figure}[ht]
\center{\includegraphics[width=0.95\linewidth]{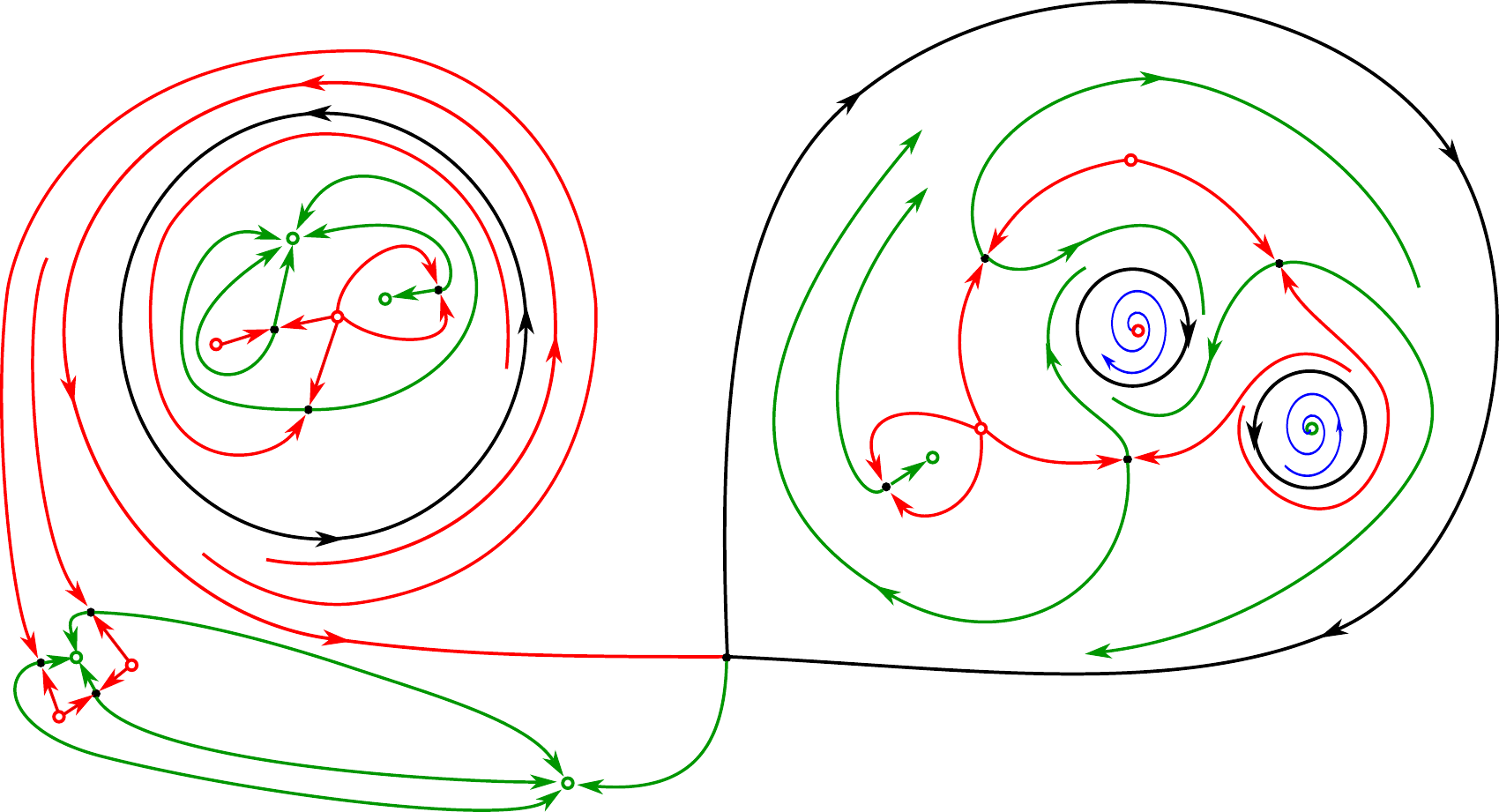}}
\caption{Separatrix diagram of the $HC_-$-flow on the two-dimensional sphere}
\label{fig:exam}
\end{figure}

We can see a diagram for a specific $HC_-$-flow on the two-dimensional sphere in Fig.~\ref{fig:exam}. It is a rough sketch of the flow where its closed trajectories, singular points, the homoclinic trajectory, stable and unstable manifolds of all the saddles are highlighted.

Its structural graph and distinguishing graphs are depicted below.

\begin{figure}[ht]
\center{\includegraphics[width=0.95\linewidth]{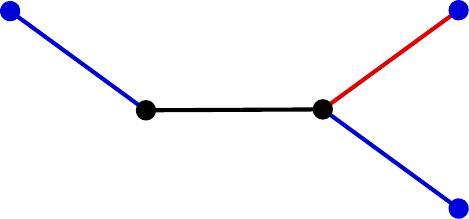}
}
\put (-295,70) {\large +}
\put (-133,70) {\large --}
\put (-360,120) {\large $1$}
\put (-360,160) {\large --}
\put (-214,110) {\large +}
\put (-68,120) {\large $2$}
\put (-68,160) {\large --}
\put (-68,30) {\large $3$}
\put (-68,60) {\large +}
\caption{Structural graph of the $HC_-$-flow}
\end{figure}

\begin{figure}[ht]
\center{\includegraphics[width=0.40\linewidth]{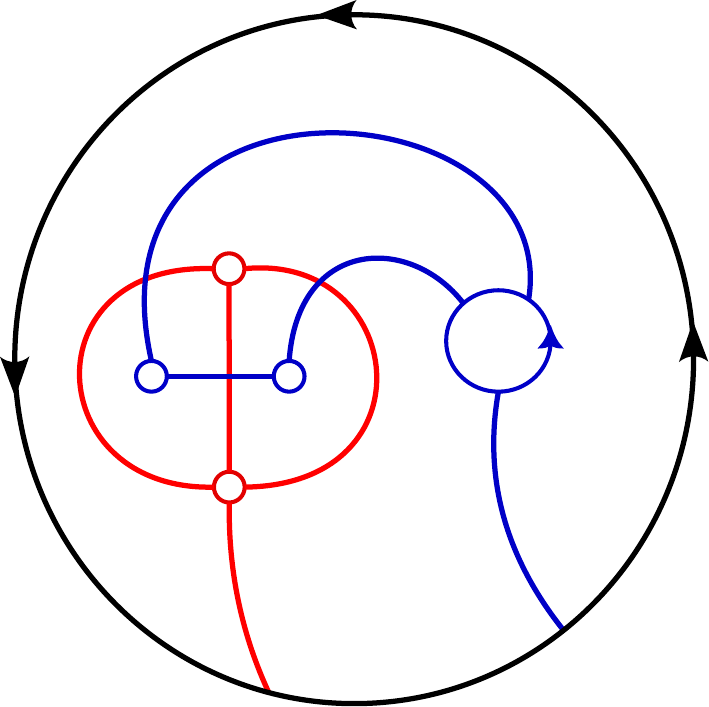}
}
\put (-5,120) {\large $0$}
\put (-40,100) {\large $1$}
\caption{Distinguishing graph of the $H_+$-atom}
\end{figure}

\begin{figure}[ht]
\center{\includegraphics[width=0.40\linewidth]{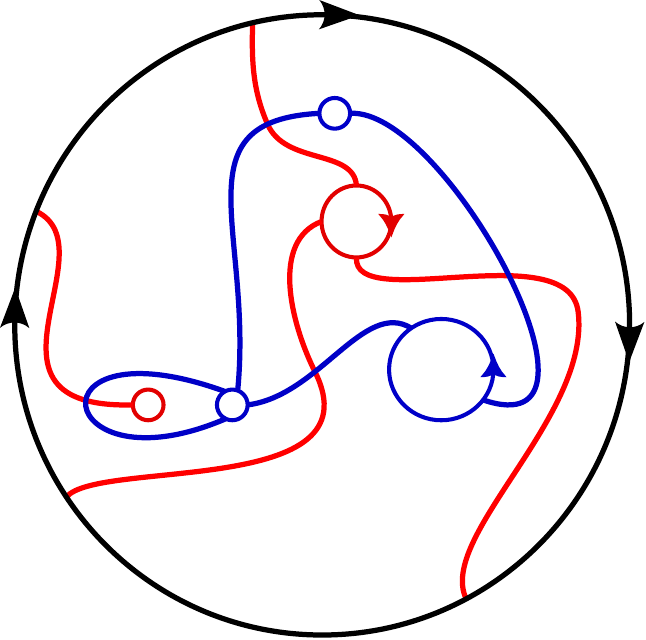}
}
\put (-5,120) {\large $0$}
\put (-68,120) {\large $2$}
\put (-45,85) {\large $3$}
\caption{Distinguishing graph of the $H_-$-atom}
\end{figure}

\begin{figure}[ht]
\center{\includegraphics[width=0.40\linewidth]{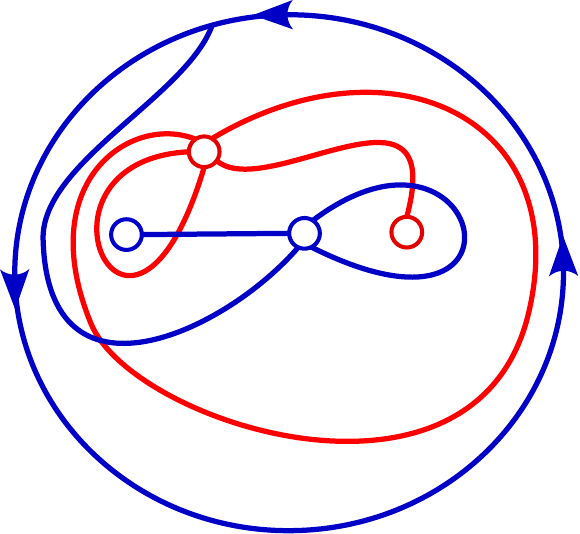}
}
\put (-5,120) {\large $1$}
\caption{Distinguishing graph of the $C$-atom adjacent to the $H_+$-atom}
\end{figure}



For the distinguishing graph describing the $H_+$-atom the blue vertex on the black edge is claimed to be a distinguishing vertex.

For the distinguishing graph describing the $H_-$-atom a blue distinguishing vertex is stated to be an upper blue vertex in oriented cycle $3$.

Finally, for the distinguishing graph describing the atom corresponding to the vertex with label $1$ in the structural graph the vertex on the blue oriented cycle is a distinguishing vertex.






The full flow code is as follows: $$S ({\color{blue}{0+_-}} {\color{red}{0-_-}} {\color{blue}{0-_+}})$$ $$D_- ({\color{blue}{\langle 1 \rangle \{123\}_{1} [24] [34]}})_{0+} ({\color{blue}{\{12\}_{3}[1344][23]}})_{0-} ({\color{blue}{\{1\}_{0+}[1223][3]}})_1$$ $$D_+ ({\color{red}{\{1234\}_{0+}\{132\}_{2}}})_{0-}$$

\section{Code as a Complete Topological Invariant}

We are going to refine the aforementioned approach to constructing flow codes, research the properties of the new codes and ways to understand whether three sequences of symbols correctly define an $HC_-$-flow on the $2$-dimensional sphere.

First, we are going to number the vertices of the structural flow code unambiguously. We leave labels $``0+``$ and $``0-``$ at their natural places. The algorithm is following.

\begin{enumerate}
  \item[(1)] Suppose that there is an observer in the point $``0+``$ standing with his back to the black edge. He numbers all the endpoints of degree $1$ from which one can find a path to the point with label $``0+``$ without passing the point with label $``0-``$ (if there are any) consecutively moving from one point to another in a counterclockwise direction starting from the least unused yet natural number and gradually increasing it by $1$.
  \item[(2)] Then, the observer shifts to the point with label $``0-``$ and numbers all the endpoints of degree $1$ from which one can find a path to the point with label $``0-``$ without passing the point with label $``0+``$ (if there are any) consecutively moving from one point to another in a counterclockwise direction starting from the least unused yet natural number and gradually increasing it by $1$.
  \item[(3)] Afterwards, the observer excludes all the numbered vertices and the edges incident on them and repeats the operation until he runs out of unlabeled vertices.
\end{enumerate}

The obtained numbering is called a \emph{canonical numbering}. There is the only numbering of closed trajectories in an $HC_-$-flow for which the function $f$ on the set of the vertices defined in the previous chapter gives us the canonical numbering. Hereafter, we will presume that the vertices of the structural graph have always the canonical code and that the oriented cycles of the distinguishing graphs are labeled in a relevant way.


\begin{definition}
An $HC_-$-flow on the two-dimensional sphere with $N$ closed trajectories, $N \in \mathbb{N}$, is called \emph{distinguished} if for each its atom some single non-homoclinic trajectory belonging to a stable manifold of some saddle point if it exists and its alpha-limit set are fixed.
\end{definition}

Thus, the corresponding fixed trajectory in each atom is called \emph{distinguished} as well.

In this way, we automatically fix in the distinguishing graph for each atom a \emph{distinguishing vertex}.

\begin{definition}
  Suppose that $M$ is the $2$-dimensional sphere. Two distinguished $HC_-$-flows $$g_i\colon M \times \mathbb{R} \to M,\,\,i = 1,2,$$ are \emph{strongly topologically equivalent} if there exists a homeomorphism $h\colon M \to M$ which maps each trajectory of the first flow onto a trajectory of the second one without reversing orientations and maps each distinguished trajectory on $M$ into a distinguished trajectory on $M$.
\end{definition}

We fix beforehand a set $$S = \{``0+_-``, ``0-_-``, ``0+_+``, ``0-_+``, ``1_+``, ``1_-``, ``2_+``, ``2_-``, \ldots\}.$$ Let us consider any fixed finite sequence of symbols from $S$ enclosed in the string $$S(\ldots).$$ We number them as $n_1, n_2, \ldots, n_m$, $m \in \mathbb{N}$. Thus we have obtained a multiset where elements may repeat, therefore each of them should be considered with its position in the structural flow code. Then we define a function $$F\colon\{n_1, n_2, \ldots, n_m\} \to \{1, 2, \ldots, m\}$$ on the set of these symbols according to the following principle --- $F(n_i) = j$ iff $j$ is the least natural number which is not met among the symbols $n_i, n_{i+1}, \ldots, n_m$ and which is greater than $F(n_k)$, $k<i$.


\begin{theorem}\label{structure}
  A structural flow code for an $HC_-$-flow on the $2$-dimensional sphere is a string $$S(\ldots)$$ where there are $n$ symbols (not necessarily different), $n \in \mathbb{N}$, from $S$ or just an empty string instead of dots provided that in the first case the following conditions are true:
   \begin{enumerate}
     \item[(1)] all the symbols are colored either red or blue;
     \item[(2)] the last symbol of the sequence contains either $``0+``$ or $``0-``$;
     \item[(3)] the sequence does not contain natural numbers greater than $n$ and less than $2$;
     \item[(4)] if one excludes symbols $``0+``$ and $``0-``$ from the string then a string will be obtained where each natural number followed by either the same number or by this number increased by $1$;
     \item[(5)] provided that one excludes all the indices from the string, if $``0+``$ or $``0-``$ in the code is followed by another but different symbol then it never appears in the code again;
     \item[(6)] the symbol containing a natural number $m \geq 2$ must not appear later than in $m-1$-th position in the code.
   \end{enumerate}
\end{theorem}

\begin{proof}
  The first statement in the theorem is obvious due to the described algorithm of constructing the structural flow code.

  As for the second statement, according to the algorithm, the last symbol deals with the only left non-black edge which is attached to the vertex with either a label $``0+``$ or a label $``0-``$. Hence, the second statement is true as well.

  One also has no difficulties in seeing a reason why the third statement is true --- if the structural flow code has $n$ symbols then it has $n$ vertices numbered from $1$ to $n$ and two vertices assigned with labels $``0+``$ and $``0-``$. Besides, there is no way $1$ can appear in the code what follows from the canonical numbering.

  The fourth, the fifth and the sixth statements are true what follows from the way of constructing a canonical numbering of the structural graph.
\end{proof}

\begin{theorem}\label{first}
  The first string in a distinguishing flow code for an $HC_-$-flow on the $2$-dimensional sphere is a string $$D_-\ldots$$ where there is a concatenation of multiple $HC$-like codes placed in ascending order of the indices near their closing brackets instead of dots with the following properties:
  \begin{enumerate}
    \item[(1)] the concatenation certainly contains an $HC$-like code with an index $0+$ and this code is the only one which contains angle brackets;
    \item[(2)] there is only one number inside angle brackets;
    \item[(3)] each symbol inside of each $HC$-like code list is colored blue;
    \item[(4)] each natural number present in an arbitrary fixed $HC$-like code is placed there twice;
    \item[(5)] in any $HC$-like code each pair of adjacent $HC$-like code lists in this sequence contains a pair of the same natural numbers placed in different code lists;
    \item[(6)] if there is an $HC$-like code $A$ with its corner index $i \in S$ which contains an $HC$-like code list in curly brackets with some symbol $j \in S$ in its right lower corner then provided that there is another $HC$-like code $B$ with its corner index $j \in S$, $B$ contains a code list in curly brackets with the symbol $i \in S$ in its right lower corner;
    \item[(7)] if $V$ is a number of the $HC$-like code lists in an arbitrary $HC$-like code, $N$ is the greatest natural number in this very code and $F$ is a number of different cycles in this code then $$V - N + F = 2.$$
  \end{enumerate}
\end{theorem}

\begin{proof}
  The first three statements are quite comprehensible.

  The fourth statement is true because each natural number present in the code corresponds to an intersection of a stable manifold of a saddle (or a non-homoclinic trajectory belonging to a stable manifold of a saddle if we speak of the homoclinic saddle) with a regular neighbourhood of a source or of a repelling cycle (or with a regular neighbourhood of the union of the homoclinic saddle and a homoclinic trajectory). These intersections become intersections of boundary cycles with thin blue edges. And the same number corresponds to multiple vertices iff they are endpoints of the same thin blue edge. Thus, each natural number present in the code can appear there only twice.

  The fifth statement follows directly from connectivity of the distinguishing graphs for the atoms which in its turn follows from connectivity of the considered surface.

  The first part of the sixth statement implies that the atom corresponding to the $HC$-like code $A$ is adjacent to the atom corresponding to the $HC$-like code $B$. So, the statement about the indices is fair considering how the flow code is constructed. 


  As for the seventh statement, let us fix the distinguishing graph for an arbitrary atom which is not the $H_+$-atom and collapse each thick blue circle (which depicts a boundary of a regular neighbourhood of a source or of a repelling cycle) and each black circle (which depicts a boundary of a regular neighbourhood of a union of the homoclinic trajectory and the homoclinic saddle) into a point. All these circles correspond to the set of the code lists in the appropriate $HS$-like code. All the present thick red circles are collapsed into points as well but afterwards they are eliminated. Then we obtain a sphere with a cellular division onto it where the obtained points are $0$-cells, the thin blue edges are $1$-cells and the interiors of the cycles are $2$-cells. As the Euler characteristics of the sphere is $2$, we obtain an equality in question.

  If one deals with the $H_+$-atom then the seventh statement remains true as the algorithm is the same but one must glue to the component of the boundary of this atom containing the homoclinic trajectory $\gamma$ the $2$-disk containing a source so that all the trajectories flow from it to $\gamma$ and consider that exactly one code list here corresponds not to a source but to $\gamma$.
\end{proof}

\begin{theorem}\label{second}
  The second string in a distinguishing flow code for an $HC_-$-flow on the $2$-dimensional sphere is a string $$D_+\ldots$$ where there is a concatenation of multiple cycle-like codes for some $HC$-like codes from the first string placed in ascending order of the indices near their closing brackets instead of dots with the following properties:
  \begin{enumerate}
    \item[(1)] each symbol inside of each cycle-like code list is colored red;
    \item[(2)] a cycle-like code with an index $``0-``$ contains a cycle-like code list with an index $``0+``$;
    \item[(3)] each cycle-like code for an $HC$-like code from the first string has an index in the right lower corner which is assigned to this very $HC$-like code in the first string;
    \item[(4)] each cycle list $A$ from a cycle-like code with an index $n \in \mathbb{N}$ represents a cycle from the corresponding $HC$-like code; among all the cycle lists representing this cycle the list $A$ starts from the least possible natural number among those not included in round brackets;
    \item[(5)] cycle-like code lists in each cycle code are placed in ascending order of corresponding cycles;
    \item[(6)] if there is a cycle-like code $A$ with its corner index $i \in S$ which contains a cycle-like code list in curly brackets with some symbol $j \in S$ in its right lower corner then provided that there is another cycle-like code $B$ with its corner index $j \in S$, $B$ contains a code list in curly brackets with the symbol $i \in S$ in its right lower corner.
  \end{enumerate}
\end{theorem}


All the statements are obviously reasonable due to the construction of the second string in the distinguishing flow code.

There also have to be some specific properties for both flow codes which have to generate a structural and distinguishing graphs for $HC_-$-flows without arising any contradictions.


\begin{theorem}\label{eq1}
   Let the sequence of symbols in the round brackets in the structural code for an $HC_-$-flow $g$ on the $2$-dimensional sphere form a multiset $\{n_1, n_2, \ldots, n_m\}$. For the first string $D_-\ldots$ in a distinguishing flow code and a structural flow code $S(\ldots)$ for $g$ next conditions are equivalent provided that the string $D_-\ldots$ contains an $HC$-like code $A$ with an index $i \in S$:
   \begin{enumerate}
     \item[(1)] the code $A$ contains a code list with curly brackets and a corner index $j \in S$;
     \item[(2)] $j$ ($i$) is an element of the multiset $\{n_1, n_2, \ldots, n_m\}$ having corner index $``+``$ or $``-``$ and in one of its positions it is colored blue and mapped into $i$ ($j$) by $F$.
   \end{enumerate}
\end{theorem}

\begin{theorem}\label{eq2}
  Let the sequence of symbols in the round brackets in the structural code for an $HC_-$-flow $g$ on the $2$-dimensional sphere form a multiset $\{n_1, n_2, \ldots, n_m\}$. For the second string $D_+\ldots$ in a distinguishing flow code and a structural flow code $S(\ldots)$ for $g$ next conditions are equivalent provided that the string $D_+\ldots$ contains a cycle-like code $A$ with an index $i \in S$ and that $i, j \in S$ so that either $i \neq ``0-``$ or $j \neq ``0+``$:
   \begin{enumerate}
     \item[(1)] the code $A$ contains a code list with curly brackets and a corner index $j \in S$;
     \item[(2)] $j$ ($i$) is an element of the multiset $\{n_1, n_2, \ldots, n_m\}$ having corner index $``+``$ or $``-``$ and in one of its positions it is colored red and mapped into $i$ ($j$) by $F$.
   \end{enumerate}
\end{theorem}


\begin{theorem}\label{blank}
  Let the sequence of symbols in the round brackets in the structural code for an $HC_-$-flow $g$ on the $2$-dimensional sphere form a multiset $\{n_1, n_2, \ldots, n_m\}$. Suppose that the first string $D_-\ldots$ in a distinguishing flow code for $g$ does not contain an $HC$-like code with its index $i \in S$. Then if $i$ is present in the structural flow code only one of three following options is possible:
  \begin{enumerate}
    \item[(1)] $i = ``0+``$ and it appears in the structural flow code only once colored blue;
    \item[(2)] $i$ depicts some natural number not greater than $m$ and it appears in the structural flow code only once; besides, another symbol having the other colour exists in the code so that $F$ maps it into $i$.
  \end{enumerate}
\end{theorem}

\begin{definition}
  The ordered set of the three strings of symbols are called an \emph{$HC$-defining triple} if they are enclosed in the following templates: $$S(\ldots),\,\,\,\,D_-\ldots,\,\,\,\,D_+\ldots.$$ The symbols placed instead of the dots in each string are called a \emph{kernel} of the appropriate string.
\end{definition}


\begin{theorem}[Realisation theorem]
  Suppose that there is an $HC$-defining triple. Furthermore, for the first string the list of conditions from Theorem~\ref{structure} is true, for the second string the list of conditions from Theorem~\ref{first} is true, and for the third string the list of conditions from Theorem~\ref{second} is true. Besides, the statements from Theorem~\ref{eq1}, Theorem~\ref{eq2} and Theorem~\ref{blank} are true as well.

  Then, these three strings form a structural and a distinguishing flow code for a specific distinguished $HC_-$-flow on the $2$-dimensional sphere which is unique with regard to strong topological equivalence of $HC_-$-flows.
\end{theorem}

\begin{proof}
  Three codes give us unambiguously the structural graph for the flow and all the distinguishing graphs for all its atoms. If code lists for some atoms are absent, it means these atoms contain no saddles and the structure of their distinguishing graphs is implied from what is written in the structural flow code. Furthermore, we can restore all the labels assigned to the oriented cycles on the graphs and tell that two oriented cycles are equivalent iff they are assigned by the same number.

  Now, for each pair of equivalent oriented cycles if they depict an attracting closed trajectory we take a system of differential equations on the plane dealing with polar coordinates:
  \begin{equation*}
\left \{
\begin{aligned}
  & \displaystyle\frac{dr}{dt} = r(1-r^2), \\
  & \displaystyle\frac{d\Theta}{dt} = \pm 1.
\end{aligned} \right.
\end{equation*}

Then, we take a regular neighbourhood of the closed trajectory placed along the curve $r = 1$. The sign in the second equation depends on which direction the movement on this trajectory is moving. If we are interested in a repelling closed trajectory, we take a bit different system of differential equations:
\begin{equation*}
\left \{
\begin{aligned}
  & \displaystyle\frac{dr}{dt} = r(r^2-1), \\
  & \displaystyle\frac{d\Theta}{dt} = \pm 1.
\end{aligned} \right.
\end{equation*}

The structural graph gives us a configuration of these trajectories relatively to each other. Also we are taking a regular neighbourhood of a union of a homoclinic trajectory and a homoclinic saddle of an appropriate dynamical system on the Euclidean plane --- it is depicted by a pair of black oriented cycles.

Afterwards, in each atom we replace a non-oriented blue boundary cycle with a $2$-disk $$B = \{(x, y) \in \mathbb{R}^2\,|\,x^2+y^2\le 1\}$$ with the specified flow generated by the vector field $$Y_1=\{x,y\}$$ on it. Each thin blue edge is replaced by a region $$W = \{(x, y) \in \mathbb{R}^2\,|\,-1\le x \le 1, -1\le y \le 1\}$$ with the specified flow generated by the vector field $$Y_2=\{-x,y\}$$ on it. Finally, each canonical $2$-cell without a red oriented cycle inside it is replaced by a $2$-disk $B$ but with the specified flow inverse to the defined one on it before. Each empty region is filled by a $2$-disk $B$ but there is another specified flow on it and it is generated by the vector field $$Y_3 = \{1,0\}$$ so that the trajectories were directed towards attracting cycles. We got a phase portrait of a smooth dynamical system on the $2$-dimensional sphere after consecutive smoothing of the vector firlds right before gluing. 
\end{proof}


\begin{theorem}[Strong topological equivalence criterion]
  Two distinguished $HC_-$-flows on the $2$-dimensional sphere are strongly topologically equivalent iff their triples of structural flow codes and parts of distinguishing flow codes are identical.
\end{theorem}

Now, let us research a phenomenon of ordinary topological equivalence of such kind of flows. We will assume that the sets $\Lambda_1, \Lambda_2, \Lambda_3, \Lambda_4, \Lambda_5$ have the same definition as they did in the previous chapter. Also, for an ordered multiset $\{a_1, a_2, \ldots, a_n\}$ we denote $\tau_{a_i}(\{a_1, a_2, \ldots, a_n\})$ its cyclic permutation starting from $a_i$ where $i$ is fixed number from $\{1, 2, \ldots, n\}$. Analogously if one has an $HC$-like code list $$H = q_1 a_1 a_2 \ldots a_n q_2$$ with non-numerical symbols $q_1$ and $q_2$ and perhaps some other non-numerical symbols whose quantity is less than $n$ it will be said without arising any contradictions that $$\tau_{a_i}(H) = q_1 L q_2$$ where $L$ forms a multiset $\tau_{a_i}(\{a_1, a_2, \ldots, a_n\})$ (as earlier, two-digit numbers must be separated with commas from each side).



Let us presume that we have a sequence of code lists $$H_1 H_2 \ldots H_n.$$ We will define a \emph{listing linear order} on the set of numerical symbols included in this list so that $a \leq b$ for two numerical symbols $a$ and $b$ from this sequence iff the first apparition of $a$ happens earlier in this sequence than the first apparition of $b$.

\begin{definition}\label{thesame}
  Let $$q_1 H_1 H_2 \ldots H_n q_2$$ be an $HC$-like code including $HC$-like code lists $H_1, H_2, \ldots, H_n$ and non-numerical symbols $q_1$ and $q_2$. Then a new $HC$-like code $$q_1 L_1 L_2 \ldots L_n q_2$$ with $HC$-like code lists $q_1$ and $q_2$ is called an \emph{equivalent $HC$-like code regarding to the symbol $p_1 \in H_i$} where $p_1$ and $H_i$ are fixed, $i = \overline{1, n}$, if it is created according to the next algorithm:
  \begin{enumerate}
    \item[(1)] first, $Q_1 = \tau_{p_1}(H_i)$, and for $s = \overline{2, n}$, $Q_s = \tau_{p_s}(H_l)$, $l = \overline{1, n}$, iff due to the listing linear order defined on the sequence of code lists $Q_1 Q_2 \ldots Q_{s-1}$ there exists the least numerical symbol which appears in this sequence only once but its second copy $p_s$ appears in $H_l$;
    \item[(2)] second, $L_1 L_2 \ldots L_n$ is formed from $Q_1 Q_2 \ldots Q_n$ through writing it down but replacing numerical symbols in the following way --- each copy of a number $p$ is replaced by a copy of a number $m$ which means that the first apparition of $p$ in $Q_1 Q_2 \ldots Q_n$ has number $m$ among all first apparition of the symbols.
  \end{enumerate}
\end{definition}

Now, let us assume that a cycle-like code list $A$ for the $HC$-like code $$H_1 H_2 \ldots H_n$$ is fixed. Suppose that the relevant cycle list is starting with some numerical symbol $q \in H_i$, where $H_i$ is fixed, $i = \overline{1, n}$. After the actions done in the items $(1)$ and $(2)$ in the Definition~\ref{thesame}, this symbol was moved into some position in some fixed $HC$-like code list $Q_j$, $j = \overline{1,n}$. Then, the third action turned this symbol into another numerical symbol $b$ belonging to $Q_j$. Thus, the cycle-like code list $B$ relevant to the cycle list starting with $b \in Q_j$ will be called \emph{relevant to the cycle-like code list $A$ through the equivalence} mentioned in Definition~\ref{thesame}.

\begin{definition}
  Suppose that $q_1 A_1 A_2 \ldots A_n q_2$ and $q_1 B_1 B_2 \ldots B_n q_2$ are cycle-like codes where $q_1$ and $q_2$ are some non-numerical symbols. Then the second cycle-like code is called \emph{relevant to the first cycle-like code through some equivalence} if for each $i = \overline{1,n}$ $A_i$ and $B_i$ are relevant through the same equivalence.
\end{definition}

It is worth indicating that if one has a distinguishing graph $\Gamma$ for any atom in an $HC_-$-flow on the $2$-dimensional sphere and has constructed the code $$H_1 H_2 \ldots H_n$$ mentioned in Definition~\ref{thesame} just as the flow code for $\Gamma$ with some distinguished vertex then the shift to the equivalent $HC$-like code from this definition is merely constructing another flow code for $\Gamma$ but with another distinguishing vertex represented in the initial flow code through the symbol $p_1 \in H_i$! Moreover, the shift to the relevant cycle-like code list through this very equivalence is just a reencoding some canonical cycle in the atom using new flow code.

Thus, it is correct to speak of two \emph{equivalent $HC$-like codes} if any of them is equivalent to the other one regarding to the specific numerical symbol in the specific $HC$-like code list.

\begin{example}
  Let us consider the $HC$-like code $$(\langle 1 \rangle \{123\}_{1} [24] [34])_{0+}.$$ The set of all $HC$-like codes equivalent to it looks like this: $$(\langle 1 \rangle \{123\}_{1} [24] [34])_{0+},\,\,(\{123\}_{1} \langle 1 \rangle  [24] [34])_{0+},\,\,(\{123\}_{1} [14] [24] \langle 3 \rangle)_{0+},$$ $$(\{123\}_{1} [14] \langle 2 \rangle [34])_{0+},\,\,([12] \{134\}_{1} [23] \langle 4 \rangle)_{0+},\,\,([12] [13] \{234\}_{1} \langle 4 \rangle)_{0+},$$ $$([12] \{134\}_{1} [24] \langle 3 \rangle)_{0+},\,\,([12] [13] \{243\}_{1} \langle 4 \rangle)_{0+}.$$
\end{example}

\begin{theorem}[Topological equivalence criterion]
  Suppose that $$g_i\colon M \times \mathbb{R} \to M,\,\,i = 1,2,$$ are two $HC_-$-flows defined on the $2$-dimensional sphere $M$ having their structural and distinguishing flow codes expressed as $HC$-defining triples $\Omega_i$, $i = 1, 2$. Then these flows are topologically equivalent iff three following statements are true:
  \begin{enumerate}
    \item[(1)] the first strings in $\Omega_1$ and $\Omega_2$ are identical;
    \item[(2)] kernels of the second strings in $\Omega_1$ and $\Omega_2$ form respectively concatenations of $HC$-like codes $L_1 L_2 \ldots L_n$ and $S_1 S_2 \ldots S_n$ so that for each $i = \overline{1,n}$ $L_i$ is equivalent to $S_i$;
    \item[(3)] kernels of the third strings in $\Omega_1$ and $\Omega_2$ form respectively concatenations of cycle-like codes $Q_1 Q_2 \ldots Q_m$ and $C_1 C_2 \ldots C_m$, $m \leq n$ where $n$ is the number indicated in the previous item and for each $j = \overline{1,m}$ $C_j$ is relevant to $Q_j$ through the equivalence mentioned in the previous item.
  \end{enumerate}
\end{theorem}

Two $HC$-defining triples whose properties described in items $(1)-(3)$ in this theorem are called \emph{equivalent}.

\section{Examples}

Let us introduce another definition.

\begin{definition}
  Two $HC_-$-flows defined on the $2$-dimensional sphere are called \emph{structurally equivalent} if they have identical structural flow codes.
\end{definition}

It is pretty obvious that from strong or ordinary topological equivalence of such flows follows their structural equivalence. Notwithstanding, their structural equivalence does not yet implies any kind of their topological equivalence.

It is possible to distinguish how many $HC_-$-flows up to their structural equivalence one can construct provided these flows have at most $6$ closed trajectories. Of course, tangible results will be obtained by means of properties of structural flow codes researched in the previous section (first of all, Theorem~\ref{structure} is needed). So, let us realize the magnitude of the problem.

The flow may have no trajectories at all, and then the kernel of its structural graph is empty. There is only $1$ such flow up to structural equivalence.

If the flow has $1$ closed trajectory, the kernel of its structural flow code may contain $``0+``$ or $``0-``$ and no natural numbers at all. This symbol has to be colored either red or blue. As a result --- $8$ flows.

If the flow has $2$ closed trajectories, the kernel of its structural flow code does not contain natural numbers distinct from $2$. Overall, we have $96$ flows.

If the flow has $N$ closed trajectories, $N \geq 3$, the kernels of its structural flow code does not contain natural numbers greater than $N$ and less than $2$. Here one may need more combinatorial methods to calculate the number of these flows.

Having $N$ symbols in the structural flow code, we are about to count the number of ways to place $``0-``$ and $``0+``$ in this code. This number is obviously calculated like this without considering so far that each symbol is supposed to have a specific colour and a specific index: $$2 \cdot \sum\limits_{k=1}^N (1+ \sum\limits_{j\geq 1, j \leq N-k} j).$$

After we placed symbols $``0+``$ and $``0-``$ in any way in some places of the kernel of the structural flow code some places remain vacant but we remember their order number in the flow code. Henceforth, we introduce a symbol $$f^{N, k}_{\{x_1, \ldots, x_j\}}$$ where $k$ and $N$ denote natural numbers, $k \leq N$, and $x_1, \ldots, x_j$ denote consecutive natural numbers not greater than $N-k$ (mind you, the set $\{x_1, \ldots, x_j\}$ may certainly be empty). It denotes a number of ways to write a complete form of the kernel of the structural flow code with $N$ symbols where $k$ last symbols and the symbols written in positions from $x_1$ to $x_j$ are filled with $``0+``$ and $``0-``$.

Thus, the number of all the ways of writing the structural flow code with $N$ symbols is equal to
$$S_N = 4^{N+1} \cdot \sum\limits_{k=1}^N \left(\sum\limits_{j\geq 0, j \leq N-k} \left(\underbrace{f^{N, k}_{\{1, \ldots, j\}} + f^{N, k}_{\{2, \ldots, j+1\}} + \ldots + f^{N, k}_{\{N-k-j+1, \ldots, N-k\}}}_{N-k-j+1\,\mathrm{summands}}\right)\right).$$

Note that for an arbitrary fixed natural number and for any $k \leq N$ $$f^{N,k}_{\{1, 2, \ldots, N-k\}} = 1.$$

The problem of direct calculating of these structural flow codes requires more advanced methods so for simplicity we consider small enough numbers of closed trajectories. We must calculate all the appropriate values of $f^{N, k}_{\{x_1, \ldots, x_j\}}$ so as to evaluate the overall number of the flows for a specific $N$.

$\textbf{N=3}$. We need to analyze the following values: $$f^{3,2}_{\varnothing},\,\,f^{3,1}_{\varnothing},\,\,f^{3,1}_{\{1\}}.$$ $f^{3,1}_{\{1\}}$ denotes a number of ways to fill the second position in a three-position structural flow code if all the other positions are already filled. The option is unique --- we can put only $3$. So, $$f^{3,1}_{\{1\}} = 1.$$ Analogously, $$f^{3,2}_{\varnothing} = 2,\,\,f^{3,1}_{\varnothing} = 2,\,\,f^{3,1}_{\{2\}} = 2.$$ Thus, $$S_3 = 256 \cdot ((1 + 1 + 2 + 2) + (1 + 2) + 1) = 2560.$$

$\textbf{N=4}$. Now we are interested in the following values: $$f^{4,1}_{\varnothing} = 4,\,\,f^{4,1}_{\{1\}} = 2,\,\,f^{4,1}_{\{2\}} = 2,\,\,f^{4,1}_{\{3\}} = 4,$$
$$f^{4,1}_{\{1,2\}} = 1,\,\,f^{4,1}_{\{2,3\}} = 3,\,\,f^{4,2}_{\varnothing} = 4,\,\,f^{4,2}_{\{1\}} = 2,\,\,f^{4,2}_{\{2\}} = 3,\,\,f^{4,3}_{\varnothing} = 3.$$

Using the aforementioned formula, we are getting the following result:
$$S_4 = 1024 \cdot (17 + 10 + 4 + 1) = 32768.$$

So, due to structural equivalence of $HC_-$-flows, you can generate $1$ such flow with no closed trajectories, $8$ flows with $1$ closed trajectory, $96$ flows with $2$ closed trajectories, $2560$ flows with $3$ closed trajectories and $32768$ flows with $4$ closed trajectories.

Now it is worth researching a number of homoclinic flows on the $2$-dimensional sphere due to topological equivalence.

Let $$g\colon M \times \mathbb{R} \to M$$ be a Morse flow on a $2$-dimensional sphere $M$. Suppose that $x_0$ is an arbitrary source of this flow and $U_{x_0}$ is its regular neighbourhood. Then, considering that $\Omega_1(g) \neq \varnothing$ trajectories which come out of $x_0$ and flowing into saddles of $g$ divide $U_{x_0}$ into \emph{sectors}. One may construct a $HC_-$-flow following these steps:
\begin{enumerate}
  \item[(1)] First, we have to construct the $H_+$-atom. We take $x_0$ and an arbitrary sector near it. It defines unambiguously a sink $s_0$ of this flow, or, more precisely, its cellular neighbourhood $U$ whose boundary is formed by trajectories belonging to stable manifolds of saddles (or the saddle if it is unique) of $g$. We transform $\mathrm{Int}(U)$ so that there is a new homoclinic saddle $u_0$ inside this neighbourhood with its homoclinic trajectory, the trajectory coming out of $x_0$ into $u_0$ and the trajectory coming out of $u_0$ into $s_0$.
  \item[(2)] Second, we construct the $H_-$-atom. We have to fill in the space inside the homoclinic trajectory. In order to do that, we should take an arbitrary Morse flow $G$ on a $2$-dimensional sphere $N$ and to transform an arbitrary sink in $G$ into a hole in the sphere which contains an attracting closed trajectory. Thus, we glue the aforementioned $H_+$-atom and the newly constructed flow on the sphere with $1$ hole.
  \item[(3)] If one intends to add some $C$-atoms, he has to do the second step for any source or sink in $g$ keeping in mind that a source must be replaced by a $2$-disk having a \emph{repelling} closed trajectory and a sink must be replaced by a $2$-disk having an \emph{attracting} closed trajectory.
\end{enumerate}

Let $E_k^m$ denote the number of choices of Morse flows on the $2$-dimensional sphere with a marked
oriented sector near a source (due to topological equivalence) and $k$ saddles and markings of $m$ sources or sinks in them. It is equivalent to calculation of the number of all possible $H_+$-atoms with $k$ non-homoclinic saddles and $m$ closed trajectories.

Let also $B_k^m$ denote the analogous number of choices of Morse flows on the $2$-dimensional sphere with a marked
source together with a given orientation at the boundary of its regular neighbourhood. It is, in its turn, equivalent to calculation of the number of all possible ways to construct $H_-$-atoms and $C$-atoms with $k$ non-homoclinic saddles and $m$ closed trajectories.

Besides, we introduce notations $E_k := E_k^0$, $B_k := B_k^0$.

The following formula holds true:
\[
E_k^m=\binom{k+2}{m}E_k^0.
\]

Let $N_k^i$ be a number of homoclinic flows on the $2$-dimensional sphere with $k$ saddles and $i$ closed trajectories.

\begin{theorem}
  The following formulas hold true: $$N_k^0(S^2) = 2 \sum\limits_{i=0}^{k-1} E_i B_{k-i-1},$$ $$N_k^1(S^2) = 4 \sum\limits_{i+j+m=k-1} \left(E_i B_j^1 B_m + E_i^1 B_j B_m \right),$$ $$N_k^2(S^2) = 8 \sum\limits_{i+j+p+q=k-1} \left(2 E_i^1 B_j^1 B_p B_q + E_i B_j^1 B_p^1 B_q + E_i B_j^2 B_p B_q + E_i^2 B_j B_p B_q \right),$$
  \[
  \begin{aligned}
  N_k^3(S^2)
  &=
  16
  \sum\limits_{i+j+p+q+r=k-1}
  \Bigl(
  3E_i^1 B_j^1 B_p^1 B_q B_r + E_i B_j^1 B_p^1 B_q^1 B_r + 2 E_i B_j^1 B_p^2 B_q B_r + \\
  &\quad + 2 E_i^2 B_j^1 B_p B_q B_r + 2 E_i^1 B_j^2 B_p B_q B_r + E_i^3 B_j B_p B_q B_r + E_i B_j^3 B_p B_q B_r \Bigr).
  \end{aligned}
  \]
\end{theorem}

After obtaining these formulas it is more convenient to go on with further calculations which help us to comprehend the number of homoclinic flows for some small numbers of saddles and closed trajectories in it. By means of the program code \cite{https://doi.org/10.5281/zenodo.22298943} we complete calculations for small numbers of $k$ and $m$ thus obtaining the following table of results:

\begin{table}[ht]
\centering
\begin{tabular}{c|r|r|r|r}
$k$ & $N_k^{0}$ & $N_k^{1}$ & $N_k^{2}$ & $N_k^{3}$ \\
\hline
1 & 2 & 8 & 24 & 64 \\
2 & 8 & 60 & 288 & 1072 \\
3 & 38 & 420 & 2752 & 12896 \\
4 & 220 & 3132 & 25600 & 142048 \\
5 & 1484 & 25488 & 246160 & 1558048 \\
6 & 11124 & 221920 & 2453952 & 17352576 \\
7 & 90098 & 2035500 & 25233504 & 196745568 \\
8 & 769108 & 19382868 & 265686960 & 2264046112 \\
9 & 6839220 & 189962696 & 2849218640 & 26355749664 \\
\end{tabular}
\end{table}

\section*{Declaration on the use of AI tools}

No generative AI tools were used in the preparation of this paper. The mathematical ideas, construction of the topological invariant, proofs of the main results, computational algorithms, numerical calculations, and figures are entirely the work of the authors. Both authors take full responsibility for the mathematical content and all results presented in the paper.

\bibliographystyle{plain}
\bibliography{arxiv}

\end{document}